\documentclass[11pt]{article} 

\usepackage{fullpage}

\usepackage{graphics} 
\usepackage{epsfig} 
\usepackage{amsmath} 
\usepackage{amssymb}  
\usepackage{tabularx,booktabs}
\usepackage{mathrsfs}
\usepackage{selectp}
\usepackage[normalem]{ulem}
\usepackage{cite}
\usepackage{authblk}
\usepackage{subcaption}

\usepackage{natbib}

\usepackage{scalerel,stackengine}
\stackMath
\newcommand\reallywidehat[1]{%
\savestack{\tmpbox}{\stretchto{%
  \scaleto{%
    \scalerel*[\widthof{\ensuremath{#1}}]{\kern-.6pt\bigwedge\kern-.6pt}%
    {\rule[-\textheight/2]{1ex}{\textheight}}
  }{\textheight}%
}{0.5ex}}%
\stackon[1pt]{#1}{\tmpbox}%
}

\usepackage{hyperref}
\hypersetup{
    bookmarks=true,         
    unicode=false,          
    pdftoolbar=true,        
    pdfmenubar=true,        
    pdffitwindow=false,     
    pdfstartview={FitH},    
    pdfnewwindow=true,      
    colorlinks=true,       
    linkcolor=blue,          
    citecolor=blue,        
    filecolor=magenta,      
    urlcolor=red,           
    breaklinks=true
}

\newcommand{\e}{\varepsilon}

\newcommand{\E}{\mathbf{E}}

\DeclareMathOperator{\argmin}{argmin}

\DeclareMathOperator{\tr}{tr}

\newtheorem{theorem}{Theorem}
\newtheorem{proposition}{Proposition}
\newtheorem{lemma}{Lemma}
\newtheorem{corollary}{Corollary}

\title{\LARGE \bf
How are policy gradient methods affected by the limits of control?
}

\author[1]{Ingvar Ziemann}
\author[2]{Anastasios Tsiamis}
\author[1]{Henrik Sandberg}
\author[3]{Nikolai Matni}
\affil[1]{Division of Decision and Control Systems,  KTH Royal Institute of Technology}
\affil[2]{Automatic Control Laboratory, ETH Zurich}
\affil[3]{Department of Electrical and Systems Engineering, University of Pennsylvania}
\date{}
\begin{document}

\maketitle

\begin{abstract}
We study stochastic policy gradient methods from the perspective of control-theoretic limitations. Our main result is that ill-conditioned linear systems in the sense of Doyle inevitably lead to noisy gradient estimates. We also give an example of a class of stable systems in which policy gradient methods suffer from the curse of dimensionality. Our results apply to both state feedback and partially observed systems.
\end{abstract}

\section{Introduction}
Reinforcement learning (RL) methods have shown great empirical success in controlling  complex dynamical systems \cite{silver2017mastering}. While these methods are promising, we have only begun to understand performance guarantees and fundamental limitations in continuous state and action problems. Providing such guarantees and understanding such limitations is crucial to deploying these methods in  safety-critical systems. In this paper, we focus on a particular class of such methods; namely, we seek to understand fundamental limitations for policy gradient methods.

Policy gradient methods are a relatively simple class of algorithms that have been recently analyzed in the context of the linear quadratic regulator (LQR), \cite{fazel2018global,tu2019gap}. The motivation for studying policy gradients in the context of LQR stems from that it serves as an analytically tractable benchmark for RL in continuous state and action spaces. For instance, by direct arguments on can show that control-theoretic parameters affect the hardness of both offline  and online learning in LQR \cite{tsiamis2021linear, ziemann2021uninformative, ziemann2022regret, Tasosunpub}. Here, we extend this line of work and show that the popular policy gradient methods degrade similarly for systems with poor controllability and observability. To be precise, we show that ill-conditioned systems lead to arbitrarily noisy stochastic gradients.

\paragraph{Problem Formulation}
We are interested in studying how policy gradient methods applied to the linear system
\begin{align}
    \label{eq:lds}
    x_{t+1}&=Ax_t +Bu_t + w_t, & x_0=0 &  &t=0,1,\dots
\end{align}
are affected by the fundamental limits of control. Above, $x_t \in \mathbb{R}^{d_x}, A \in \mathbb{R}^{d_x\times d_x}$, $u_t \in \mathbb{R}^{d_u}$ and $w_t \in \mathbb{R}^{d_x}$ is an i.i.d. mean zero sequence of Gaussian noise with covariance matrix $\Sigma_W \in \mathbb{R}^{d_x\times d_x}$. 

The learning task is to minimize
\begin{align}
\label{eq:costdef}
    J_S(K) \triangleq \limsup_{T\to\infty} \frac{1}{T} \sum_{t=0}^{T-1} \E_{K,S}  \left[x_t^\top Qx_t + u_t^\top R u_t\right]
\end{align}
subject to the dynamics (\ref{eq:lds}) without access to the model parameter $S= (A,B)$. In equation (\ref{eq:costdef}), $\E_{K,S}$ denotes expectation under the control law $K$ with dynamics $S$. 

In this work, we relate the efficiency of stochastic policy gradient methods to certain control-theoretic parameters. Namely, we analyze algorithms of the form
\begin{align*}
    \reallywidehat{K} \leftarrow \reallywidehat{K} -\alpha \reallywidehat{ \nabla_K J(K;S)}
\end{align*}
for some learning rate $\alpha\in \mathbb{R}_+$, and where at each iteration $\reallywidehat{ \nabla_K J(K;S)}$ is estimated using data from the system (\ref{eq:lds}). Such algorithms have been shown to converge for LQR by \cite{fazel2018global}. The purpose of this work is to demonstrate that any estimate $\reallywidehat{ \nabla_K J(K;S)}$ from an (arbitrarily) ill-conditioned system (\ref{eq:lds}) is (arbitrarily) noisy.

To make this statement rigorous we need to model the statistical information available to the learner. Here, we model this as follows: the learner is given access to $N \in \mathbb{N}$ experiments $(x_{0,n},\dots,x_{T,n}), n \in [N]$ of length $T \in \mathbb{N}$ and a total input budget of $\beta NT$ with $\beta \in \mathbb{R}_+$. More precisely, the learner is allowed to freely choose $u_{t,n}$ as a function of past observations $(x_{0,n},\dots,x_{t,n})$ and past trajectories $(x_{0,m},\dots,x_{T,m}),m<n$ and possible auxiliary randomization, while being constrained to a total budget
\begin{align}
\label{eq:budgeteq}
   \sum_{n=1}^N \sum_{t=0}^{T-1} \E u_{t,n}^\top u_{t,n} \leq \beta N T.
\end{align}
This formulation allows both open- and closed-loop experiments but normalizes the average input energy to $\beta$. 

\subsection{Related Work}

The first proof that policy gradient methods converge for LQR is given in \cite{fazel2018global}, which provides nonasymptotic guarantees that are polynomial in relevant problem parameters. Convergence guarantees for more general MDPs and other versions of LQR are given in \cite{zhang2020global, gravell2019learning, gravell2020learning, zhang2020policy,yaghmaie2022linear}. Extensions to partially observed systems are considered in \cite{tang2021analysis,mohammadi2021lack, zheng2021sample}. A popular alternative approach to policy gradients for LQR is based on certainty equivalence \cite{dean2020sample}. 

Most closely related to our work are \cite{preiss2019analyzing, tu2019gap, venkataraman2019recovering}. While we give lower bounds valid for any estimator in this work, \cite{preiss2019analyzing} analyzes the variance of a particular gradient estimator known as REINFORCE. Similarly, \cite{tu2019gap} provides algorithm specific lower bounds which demonstrate, among other things, that if $R=0$ and $B$ is invertible, then learning fails as $\|B\|\to 0$. Further, \cite{tu2019gap} gives a generic performance lower bound for offline methods which however does not scale with relevant system-theoretic quantities.

Our work also relates to \cite{tsiamis2021linear, ziemann2021uninformative, ziemann2022regret, Tasosunpub}, which also study fundamental limits in learning-enabled control. From a broader perspective, the present work fits into a line of work that strives to ascertain the interplay between control-theoretic performance, stability and robustness notions, and learning \cite{bernat2020driver, boffi2021regret, perdomo2021stabilizing, tu2021sample,ziemann2022single}. While we are mostly interested in offline methods in this work, analyses of online LQR can be found in the literature, see \cite{abbasi2019model, simchowitz2020naive,cassel2021online, ziemann2022regret} and the references therein.

\subsection{Contribution}
We show that policy gradient methods are very much affected by the limits of control. Our main result ( Theorem~\ref{thm:mainthm}) demonstrates that state feedback systems operating near marginal stability suffer from noisy gradients. This happens for instance if the system has poorly controllable unstable modes. We also provide an analogue of this result for partially observed systems (Theorem~\ref{thm:mainthmpo}), which we use to show that systems with bad (small) Markov parameters also lead to noisy gradients. Compared to previous literature on this topic \cite{tu2019gap, venkataraman2019recovering,preiss2019analyzing}, our results provide a more fine-grained theoretical understanding of when and why gradient methods applied to dynamical systems fail.

\subsection{Preliminaries} A matrix $A$ is stable if $\rho(A) <1$. A matrix $K$ is stabilizing for the system $S=(A,B)$ if $A+BK$ is a stable matrix. If $K$ is stabilizing for the system $S=(A,B)$, the closed-loop controllability gramian
\begin{align}
\label{eq:gammadef}
    \Gamma_{K,S} \triangleq  \sum_{t=0}^\infty (A+BK)^t\Sigma_W(A+BK)^{t,\top}
\end{align}
is well-defined. The set of all systems $S$ for which there exists a stabilizing $K$ is denoted by $\mathcal{S}=\mathcal{S}_{d_x,d_u}$, which is an open subset of $\mathbb{R}^{d_x \times d_x + d_x\times d_u}$ in norm topology. Further, we define $K_\star(S)$ as (any element of) $    K_\star(S) \in \argmin_K J_S(K) $. Moreover, we denote the matrix operator norm (induced $l^2(\mathbb{R}^d) \to l^2(\mathbb{R}^d)$) by $\| \cdot \|_{\mathsf{op}}$.

We also require the following information-theoretic quantities. We define the Kullback-Leibler divergence between two probability measures $\mathbf{P}$ and $\mathbf{Q}$ as $d_{\mathsf{KL}}(\mathbf{P},\mathbf{Q})\triangleq\int \frac{d\mathbf{P}}{d\mathbf{Q}} d\mathbf{P} $ and the total variation distance as $d_{\mathsf{TV}}(\mathbf{P},\mathbf{Q})\triangleq\int |d\mathbf{P}- d\mathbf{Q}|$. When $\mathbf{P}$ and $\mathbf{Q}$ correspond to induced probability measures from two systems $S_1$ and $S_2$ we abuse notation and write $d_{\mathsf{KL}}(S_1,S_2) = d_{\mathsf{KL}}(\mathbf{P},\mathbf{Q})$ and $d_{\mathsf{TV}}(S_1,S_2) = d_{\mathsf{TV}}(\mathbf{P},\mathbf{Q})$ for divergences between the corresponding parametric families.

It will be convenient to introduce the shorthand $a_t \lesssim b_t$ if there exists a universal constant $C$ such that $a_t \leq C b_t$ for every $t \geq t_0$ and some $t_0 \in \mathbb{N}$. If $a_t \lesssim b_t$ and $b_t \lesssim a_t$ we write $a_t \asymp b_t$. For an integer $N$, we also define $[N] \triangleq \{1,\dots,N\}$.
\paragraph{Policy Gradients}
We begin by recalling a standard characterization of the LQR cost (\ref{eq:costdef}) for a linear controller $K$. A version of the following lemma can also be found in for instance \cite{fazel2018global}.
\begin{lemma}
\label{lem:lqrcostlem}
If $K$ is stabilizing for $S=(A,B)$, the LQR cost can be written as
\begin{align*}
    J(K;S) = \tr P_K \Sigma_W
\end{align*}
where $P_K$ satisfies the Lyapunov equation
\begin{align}
\label{eq:lyapk}
    P_K = Q+K^\top RK + (A+BK)^\top P_K (A+BK).
\end{align}
\end{lemma}

Lemma~\ref{lem:lqrcostlem} allows us to conveniently characterize the policy gradient $\nabla_K J(K;S)$.

\begin{lemma}
\label{lem:pgexplicitlemma}
Let $K$ be stabilizing for $S=(A,B)$. The policy gradient $\nabla J(K;S)$ can be written as
\begin{align}
  \label{eq:pgexact}
    \nabla_{ K} J(K;S)= 2\left((R+B^\top P_K B)K+B^\top P_K A \right) \Gamma_{K,S}
\end{align}
where $P_K$ satisfies the Lyapunov equation (\ref{eq:lyapk}) and where $\Gamma_{K,S}$ is given by definition~(\ref{eq:gammadef}).
\end{lemma}

Combining Lemmas~\ref{lem:lqrcostlem} and \ref{lem:pgexplicitlemma} we see that we are almost in the same setting as studied in \cite{fazel2018global}. The difference is mainly in how information is acquired, since each sample from the system (\ref{eq:lds}) is noisy and conditionally Gaussian. Compared to the noise-free random initial condition setting considered in \cite{fazel2018global}, this simplifies the analysis of the variance of the gradient estimates since we later rely on the closed form of the KL divergence for Gaussians of different means.

\section{Policy Gradient Estimation Lower Bounds}
\label{sec:mainsec}

Let us begin our study of stochastic gradient methods by the observation that $\nabla_K J(K;S)$ diverges if $K$ does not stabilize the system (\ref{eq:lds}). Consider for instance the following two systems
\begin{equation}
\label{eq:twoscalars}
    \begin{cases}
 S_1:   x_{t+1} = a x_t + bu_t +w_t,\\
 S_2:   x_{t+1} = ax_t- bu_t+w_t.
    \end{cases}
\end{equation}
If $|a|>1$ there exists no linear feedback controller which stabilizes both $S_1$ and $S_2$ of equation (\ref{eq:twoscalars}). Hence, any policy gradient which is finite for the first system will be infinite for the second system and vice versa. Combining this observation with the two point method (Lemma~\ref{thm:lecam}) leads to the following conclusion.

\begin{proposition}
\label{prop:blowupbound}
For any $k \in \mathbb{R}$ the global minimax complexity of estimating the policy gradient at $k$ is infinite:
\begin{equation}
\label{eq:blowupbound}
    \inf_{\reallywidehat{\nabla J} }\sup_{(a,b) \in \mathbb{R}^2} \E_{(a,b)} \left| \frac{d}{dk} J(k;(a,b))-\reallywidehat{\nabla J} \right| = \infty
\end{equation} 
where the infimum is taken over all measurable functions of the data $(x_{0,n},u_{0,n},\dots,u_{T-1,n},x_{T,n}), n \in [N]$.
\end{proposition}

Proposition~\ref{prop:blowupbound} shows that the global minimax complexity of estimating gradients is infinite. While this shows that estimating gradients can be hard, it does not reveal how this hardness depends on control theoretic parameters. 

\subsection{Local Minimax Complexities}

In order to understand what properties of a particular system makes learning to control hard, we need to consider  local complexity measures. Here, we investigate the $(d,\e)$-local minimax complexity of estimating gradients. We define this as
\begin{equation}
\label{eq:minimaxdef}
    \mathfrak{M}_d(\e;S,K) \triangleq \inf_{\reallywidehat{\nabla J} }\sup_{S': d(S,S')\leq \e} \E_{S'} \left\| \nabla_K J(K;S')-\reallywidehat{\nabla J} \right\|_{\mathsf{op}}
\end{equation}
for some metric $d$ on the set of stabilizable systems $\mathcal{S}$ and where the infimum is taken over all measurable functions of the data $(x_{0,n},u_{0,n},\dots,u_{T-1,n},x_{T,n}), n \in [N]$. This captures a more instance-specific notion of how hard it is to estimate gradients. Roughly, this complexity measure corresponds to requiring algorithms to performing well not just on a nominal system $S$, but also on small $\e$-perturbations of that system. 

Note further that the definition (\ref{eq:minimaxdef}) still leaves open the question at which $K$ to measure the complexity of estimating gradients. By equation (\ref{eq:pgexact}) and Proposition~\ref{prop:blowupbound} we know that $\nabla_K J(K)$ can be arbitrarily large when evaluated far from a stationary point. As this rather reflects poor initialization than fundamental control-theoretic hardness, we instead seek to lower bound $\mathfrak{M}_d(\e;S,K)$ near $K_\star(S)$. Arguably, any successful policy gradient algorithm should eventually find itself near $K_\star(S)$. Thus, we provide lower bounds on the gradient estimation error in the vicinity of the stationary point $K_\star(S)$. We denote the associated local complexity by $\mathfrak{M}_d(\e;S)=\mathfrak{M}_d(\e;S,K_\star(S))$.

\paragraph{Constructing Hard Instances}

To simplify the evaluation of the local minimax complexity (\ref{eq:minimaxdef}), we  mainly consider the construction below. Fix a nominal instance $S_1=(A,B)$ of system (\ref{eq:lds}) with optimal control law $K_\star = K_\star(S_1)$; then the perturbation
\begin{align}
\label{eq:ldsconfuse}
    S_2:\quad  x_{t+1}&=A'x_{t}+B'u_t +w_t
\end{align}
is tractable to evaluate. Here, $A' = A - \Delta K_\star$ and $B' = B+\Delta$ for some $\Delta\in \mathbb{R}^{d_x\times d_u}$. This perturbation is convenient since $A'+B'K_\star = A+BK_\star$ for any $\Delta$ and has previously been used in \cite{simchowitz2020naive,ziemann2022regret}. In particular, the system quantities $P_{K_\star,S}$ and $\Gamma_{K_\star,S}$ are invariant as we vary $\Delta$. Combining this observation with the optimality of  $K_\star=K_\star(S_1)$ for system $S_1$, yields the following simple expression for the gradient of system (\ref{eq:ldsconfuse}).

\begin{lemma}
\label{lem:perturbpg}
The policy gradient for $S_2 =(A',B')$ given by system (\ref{eq:ldsconfuse}) at $K_\star =K_\star(S_1)$ is given by
\begin{align*}
   \nabla_K J(K_\star;S_2) = 2\Delta^\top P_{K_\star,S_1}(A+BK_\star)  \Gamma_{K_\star,S_1}
\end{align*}
\end{lemma}

\textit{Proof of Lemma~\ref{lem:perturbpg}.}
By Lemma~\ref{lem:pgexplicitlemma} the policy gradient is given by 
\begin{equation*}
  \nabla_K J(K_\star;S_2) =  2\left(RK_\star +(B+\Delta)^\top P_{K_\star,S_1}(A+BK_\star) \right) \Gamma_{K_\star,S_1}
\end{equation*}
where we used that $A+BK_\star=A'+B'K_\star$. On the other hand 
\begin{align*}
     2\left(RK_\star +B^\top P_{K_\star,S_1}(A+BK_\star) \right) \Gamma_{K_\star,S_1} =0
\end{align*}
by optimality of $K_\star$ to $S_1$. The result follows. \hfill $\blacksquare$

By combining Lemma~\ref{lem:perturbpg} with Le Cam's two point method \cite{lecam1973convergence} (provided in the appendix as Lemma~\ref{thm:lecam}) we obtain a generic estimation lower bound for policy gradients evaluated in the vicinity of the optimum $K_\star$.

\begin{theorem}
\label{thm:mainthm}
Consider two systems $S_1 =(A,B)$ and $S_2(\Delta)=(A',B')$ with $A'=A-\Delta K_\star$ and $B'=B+\Delta$. Let  $\mathfrak{M}_d(\e; S_1)  = \mathfrak{M}_d(\e; S_1,K_\star(S_1)) $ and $K_\star=K_\star(S_1)$,  then 
\begin{equation}
\begin{aligned}
\label{eq:mainthmeq}
    \mathfrak{M}_d(\e; S_1)    
    \geq \sup_{d(S_1,S_2(\Delta))\leq \e} \left\| \Delta^\top P_{K_\star,S_1}(A+BK_\star)  \Gamma_{K_\star,S_1}\right\|_{\mathsf{op}}
    \times \left(1-\sqrt{\frac{1}{2} d_{\mathsf{KL}}(S_1,S_2(\Delta))} \right).
\end{aligned}
\end{equation}
\end{theorem}

In other words, the local complexity of estimating gradients can be lower bounded by the maximum of $\left\| \Delta^\top P_{K_\star,S_1}(A+BK)  \Gamma_{K_\star,S_1}\right\|_{\mathsf{op}}$, optimized over $\Delta$ and subject to this leading to small differences in the output of systems $S_1$ and $S_2$. 

\textit{Proof of Theorem~\ref{thm:mainthm}.}
Define the loss function $L(\mathsf{dec},S) \triangleq  \left\| \nabla_K J(K;S)-\mathsf{dec} \right\|_{\mathsf{op}}$, where the decision $\mathsf{dec}$ is a placeholder variable for the gradient estimate. For any two systems   $S_1$ and $S_2$ we have that
\begin{multline*}
    L(\mathsf{dec},S_1) +L(\mathsf{dec},S_2) \\
    = \left\| \nabla_K J(K_\star;S_1)-\mathsf{dec} \right\|_{\mathsf{op}}+ \left\| \nabla_K J(K_\star;S_2)-\mathsf{dec} \right\|_{\mathsf{op}}\\
    \geq \left\| \nabla_K J(K_\star;S_1)- \nabla_K J(K_\star;S_2) \right\|_{\mathsf{op}}
\end{multline*}
by the triangle inequality. Invoking Lemma~\ref{lem:perturbpg} we thus see that for the choice $S_1 = (A,B)$ and $S_2=(A',B')$ with $A' = A-\Delta K_\star$ and $B' = B+\Delta$ we have
\begin{align}
\label{eq:Llb}
    L(\mathsf{dec},S_1) +L(\mathsf{dec},S_2)\geq \left\| 2\Delta^\top P_{K_\star,S_1}(A+BK_\star)  \Gamma_{K_\star,S_1}\right\|_{\mathsf{op}}
\end{align}
for any $\Delta \in \mathbb{R}^{d_x\times d_u}$. Combining equation (\ref{eq:Llb}) with Lemma~\ref{thm:lecam} it follows that
\begin{equation*}
\begin{aligned}
    \mathfrak{M}(\e, S_1) &\geq \left\| \Delta^\top P_{K_\star,S_1}(A+BK_\star)  \Gamma_{K_\star,S_1}\right\|_{\mathsf{op}} \left(1-d_{\mathsf{TV}}(S_1,S_2) \right)\\
    &\geq \left\| \Delta^\top P_{K_\star,S_1}(A+BK_\star)  \Gamma_{K_\star,S_1}\right\|_{\mathsf{op}}\\
    &\times \left(1-\sqrt{\frac{1}{2} d_{\mathsf{KL}}(S_1,S_2)} \right)
\end{aligned}
\end{equation*}
where the second inequality is an application of Pinsker's inequality. \hfill $\blacksquare$

At this point, we note that the right hand side of inequality (\ref{eq:mainthmeq}) is large for systems operating near marginal stability. When $A+BK_{\star}\to 1$ both $P_{K_{\star},S}$ and $\Gamma_{K_{\star},S}$ tend to infinity. To better understand the practical implications of this, we now turn to interpreting Theorem~\ref{thm:mainthm} by instantiating it for three special cases: scalar systems, over-actuated systems and integrator-like systems.

\subsection{Consequences of Theorem~\ref{thm:mainthm}}
\label{subsec:cons}
\paragraph{Scalar Systems}

The bound in Theorem~\ref{thm:mainthm} is agnostic to the experiment used to generate the dataset, which is simply reflected in the quantity $\left(1-\sqrt{\frac{1}{2} d_{\mathsf{KL}}(S_1,S_2(\Delta))} \right)$. Let us interpret Theorem~\ref{thm:mainthm} by a simple scalar example. To this end, consider the system
\begin{align}
\label{eq:sysscalar}
    s_1: x_{t+1}=ax_t+bu_t+w_t
\end{align}
with $a,b \in \mathbb{R}$, which is open-loop unstable $|a|>1$. Let $s_2$ be given by the perturbation $s_2=(a-\Delta k_\star, b+\Delta)$ with $\Delta \in \mathbb{R}$. The divergence $d_{\mathsf{KL}}(s_1,s_2)$ satisfies
\begin{equation}
\begin{aligned}
\label{eq:scalarkl}
    d_{\mathsf{KL}}(s_1,s_2) &=\sum_{n=1}^N \sum_{t=0}^{T-1}\E_{s_1} \frac{1}{2}( k_\star \Delta x_t+ \Delta u_t)^2 & &\textnormal{(Lemma~\ref{lem:KLcalc})}\\
    &\leq   \Delta^2 \sum_{n=1}^N \sum_{t=0}^{T-1}\E_{s_1} \left(  u_t^2 + k_\star x_t^2 \right) \\
    &\leq \Delta^2 NT(\Gamma_{k_\star,s_1}+\beta) & & \textnormal{(by (\ref{eq:budgeteq}))} \\
    &\leq \frac{1}{2}.
\end{aligned}
\end{equation}
if $\Delta^2 = \frac{1}{2 NT (\Gamma_{k_\star,s_1}+\beta)}$. Plugging inequality (\ref{eq:scalarkl}) into inequality (\ref{eq:mainthmeq}), we conclude that 
\begin{multline}
\label{eq:scalarsf}
    \mathfrak{M}_{d_\infty}\left( \e_{N,T},s_1 \right)\\
    \gtrsim  \frac{1}{\sqrt{N T(\beta + \Gamma_{k_\star,s_1})}}\left|  P_{k_\star,s_1}(a+bk_\star)  \Gamma_{k_\star,s_1} \right|
\end{multline}
with $\e_{N,T} \asymp \frac{1}{\sqrt{N T(\beta + \Gamma_{k_\star,s_1})}}$. In particular as $|b| \to 0$, one may verify that the right hand side of the expression (\ref{eq:scalarsf}) tends to infinity. In other words, as controllability (of unstable modes) is lost, policy gradients become arbitrarily noisy. This is verified via simulations in the appendix (Figure~\ref{fig:plots}) using both a least squares certainty equivalent approach and a $0$-th order method (see \cite[Algorithm 1]{fazel2018global}).

\paragraph{Multivariate Systems}

If we  assume that $K$ has a left nullspace, the bound in Theorem~\ref{thm:mainthm} becomes tractable to evaluate since we are free to select $\Delta$ such that $\Delta K =0$, which simplifies some calculations. Intuitively, the these instances are hard to distinguish between because controllers with left nullspaces lead to identifiability issues regarding the $B$-matrix \cite{ziemann2021uninformative}.

\begin{corollary}
\label{corr:Bcorr}
For any $\Delta$ such that $\Delta K_\star =0$ and $\|\Delta\|_{\mathsf{op}} \leq 1$ we have that
\begin{equation*}
    \begin{aligned}
    \mathfrak{M}_{d_\infty}\left( \e_{N,T},S_1 \right) \gtrsim  \frac{1}{\sqrt{\beta N T}} \left\| \Delta^\top P_{K_\star,S_1}(A+BK_\star)  \Gamma_{K_\star,S_1}\right\|_{\mathsf{op}}
    \end{aligned}
\end{equation*}
for any $\e_{N,T} \gtrsim 1/\sqrt{NT}$ and where $d_\infty(S_1,S_2) = \max (\|A-A'\|_{\mathsf{op}}, \|B-B'\|_{\mathsf{op}})$.
\end{corollary}

\textit{Proof of Corollary~\ref{corr:Bcorr}.}
Fix $\e>0$. By Lemma~\ref{lem:KLcalc} the two systems $S_1=(A,B)$ and $S_2(N,T)=(A'_{N,T},B'_{N,T})$ with $ A'=A-\frac{\e}{\sqrt{\beta NT}}\Delta K_\star =A$ and $B'=B+ \frac{\e}{\sqrt{\beta NT}}\Delta$ satisfy $d_{\mathsf{KL}}(S_1,S_2(N,T)) =O(1)$. The result follows by Theorem~\ref{thm:mainthm}. \hfill $\blacksquare$

In other words, the complexity of estimating gradients can be asymptotically lower bounded at the central limit theorem scale $\sqrt{NT}$ by the part of $P_{K_\star,S_1}(A+BK_\star)  \Gamma_{K_\star,S_1}$ that cannot be identified by closed loop experiments using $K_\star$. It so happens that this complexity measure is very similar to that dictating regret lower bounds in adaptive LQR \cite{ziemann2021uninformative, ziemann2022regret}. In the sequel, we  exploit this to show that the gradient variance can grow exponentially with the system dimension in the worse case by leveraging certain Riccati calculations due to \cite{Tasosunpub}.

\paragraph{Policy Gradients and the Curse of Dimensionality}
Let us now show that variance of policy gradient estimates can suffer from exponential complexity in the dimension. The proof of this fact relies on a construction due to \cite{Tasosunpub}. Namely, we consider a system consisting of two decoupled subsystems $S_1=(A,B)$ of the form:
\begin{equation}\label{CTRL_eq:REG_difficult_example_stable}
\begin{aligned}
 x_{t+1}&=\underbrace{\begin{bmatrix}
0&0&0&\dots&0&0\\ 0&\rho&2&&0&0\\\vdots& &&\ddots&&\vdots\\& &&\ddots&&0\\0&0&0& &\rho&2\\0&0&0& \dots&0&\rho\end{bmatrix}}_{=A}
x_t+\underbrace{\begin{bmatrix} 1&0\\0&0\\\vdots & \vdots \\ 0&0\\0&1\end{bmatrix}}_{=B}u_t+w_{t}
\end{aligned}
\end{equation}
with $\rho \in (0,1)$, $Q=I_{d_x}$ and $R=I_2$. We also define the subsystem
\begin{equation}\label{CTRL_eq:REG_difficult_example_stable_subsystem}
A_0=\begin{bmatrix}\rho&2&0&\dots&0&0\\0&\rho&2&&0&0\\\vdots& &&\ddots&&\vdots\\& &&\ddots&&0\\0&0&0& &\rho&2\\0&0&0&\dots &0&\rho\end{bmatrix},\, B_0=\begin{bmatrix}0\\0\\\vdots\\0\\1\end{bmatrix}
\end{equation}
with $Q_0=I_{d_x-1}$ and $R_0=1$ and where $A_0\in\mathbb{R}^{(d_x-1)\times (d_x-1)}$ and $B_0\in \mathbb{R}^{d_x-1}$. Note that $A_0$ is a stable matrix since $|\rho|<1$. In the notation of Theorem~\ref{thm:mainthm}, we let $\Delta=\begin{bmatrix} 0&0\\\Delta_1&0\end{bmatrix}$, so that $S_2$ consists of two weakly coupled subsystems, with coupling induced by $\Delta_1$ (recall $S_2=(A',B')=(A-\Delta K, B+\Delta)$).

Denote further by $P_{0,\star}$ the solution to the Lyapunov equation (\ref{eq:lyapk}) for the subsystem (\ref{CTRL_eq:REG_difficult_example_stable_subsystem}) with $K_0=K_{\star,0}$. Note also that $P_{0,\star}$ satisfies the discrete algebraic Riccati equation for the tuple $(A_0,B_0,Q_0,R_0)$. With these preliminaries established, we now recall the following two lemmas from \cite[Appendix E]{Tasosunpub}.

\begin{lemma}\label{CTRL_lem_app:stable_aux1}
We have:
\[
\|\Delta_1^\top P_0 (A_0+B_0K_{\star,0})\|_{\mathsf{op}} \ge \left(\frac{1}{2}+o(1)\right)(B_0'P_0B_0+R_0),
\]
where the term $o(1)$ tends to $0$ as $d_x$ tends to infinity. 
\end{lemma}

\begin{lemma}[Riccati matrix can grow exponentially]\label{CTRL_lem:exponential_riccati_stable}
For system~\eqref{CTRL_eq:REG_difficult_example_stable_subsystem} we have:
\[
B_0^\top P_{0,\star}B_0+R_0\ge 2^{2d_x-4}+1.
\]
\end{lemma}

Combining Corollary~\ref{corr:Bcorr} with Lemmas~\ref{CTRL_lem_app:stable_aux1} and \ref{CTRL_lem:exponential_riccati_stable} we arrive at the following conclusion:

\begin{proposition}
For the system $S$ given in equation (\ref{CTRL_eq:REG_difficult_example_stable}) we have that
\begin{equation}
    \begin{aligned}
    \mathfrak{M}_{d_\infty}\left( \e_T,S \right) \gtrsim  \frac{4^{d_x}}{\sqrt{\beta N T}} 
    \end{aligned}
\end{equation}
for $d_x$ and $NT$ sufficiently large.
\end{proposition}
In other words, there are classes of stable systems for which the policy gradient suffers from exponential complexity in the state dimension.

\section{Extension to Partially Observed Systems}
We now demonstrate that our lower bound approach extends to partially observed systems of the form
\begin{equation}
\label{eq:polds}
    \begin{aligned}
    x_{t+1}&=Ax_t + Bu_t + w_t, & x_0=0 &  &t=0,1,\dots\\
    y_t &= Cx_t +v_t
    \end{aligned}
\end{equation}
in which $A$ and $B$ are as in system (\ref{eq:lds}), $C \in \mathbb{R}^{d_y\times d_x}$ and both $w_t$ and $v_t$ are i.i.d. normal with mean zero and covariance $\Sigma_W,\Sigma_V$. We denote partially observed systems of the form (\ref{eq:polds}) by $G=(A,B,C)$. For system (\ref{eq:polds}) one typically seeks to learn dynamic controllers  of the form (see e.g. \cite{tang2021analysis}) 
\begin{equation}
\label{eq:poobscont}
    \begin{aligned}
    \xi_{t+1} &= A_{\mathsf{dyn}} \xi_t + B_{\mathsf{dyn}} y_t, & \xi_0=0 &  &t=0,1,\dots\\
    u_{t} &= K\xi_t
    \end{aligned}
\end{equation}
parametrized by the linear system $K_{\mathsf{dyn}}=(A_{\mathsf{dyn}}, B_{\mathsf{dyn}}, K)$. The objective, as before is to minimize the cost 
\begin{align}
\label{eq:costdefpo}
    J_G(K_{\mathsf{dyn}}) \triangleq \limsup_{T\to\infty} \frac{1}{T} \sum_{t=0}^{T-1} \E_{K_{\mathsf{dyn}},G}  \left[x_t^\top Qx_t + u_t^\top R u_t\right]
\end{align}
but this time subject to process  and controller dynamics (\ref{eq:polds})-(\ref{eq:poobscont}). 

\paragraph{Fully Observed Reformulation}
To establish a hardness result, it suffices to focus on the difficulty of estimating gradients with respect to the output matrix of the controller (\ref{eq:poobscont}), $K$. We will exploit this by reducing the system  (\ref{eq:polds})-(\ref{eq:poobscont}) when evaluated near the optimum of $J_S$ (\ref{eq:costdefpo}) to a fully observed system. Namely, at the optimum $K_{\mathsf{dyn},\star} = \argmin J_S(K_{\mathsf{dyn}})$ the dynamics of $\xi_t$ in equation (\ref{eq:poobscont}) are given by the Kalman filter $\xi_t=\hat x_t$, allowing us to write
\begin{equation}
\label{eq:filtereq}
\begin{aligned}
    \hat x_{t+1}& = A\hat x_t + B u_t +\nu_t\\
    \nu_t& \sim N(0,\Sigma_{t})
\end{aligned}    
\end{equation}
which has the same input-output behavior as system (\ref{eq:polds}) and where the sequence of innovations $\{\nu_t\}$ is independent. More precisely, the covariance $\Sigma_t$ of $\nu_t$ is given by
\begin{align}
    \label{eq:kalmancov}
    \Sigma_t = L_t ( C F_{k|k-1}C^\top +\Sigma_V)L_t^\top
\end{align}
where $F_{t|t-1}$ satisfies the filter Riccati recursion (see e.g. \cite{soderstrom2002discrete})
\begin{equation}
    \label{eq:fric}
    F_{t+1|t} =\Sigma_W+ AF_{t|t-1}A^\top  
    - F_{t|t-1}C^\top (C F_{t|t-1} C^\top +\Sigma_V)^{-1} C F_{t|t-1}
\end{equation}
and the filter gain $L_t$ is given by
\begin{align}
\label{eq:fgain}
    L_t &= F_{t|t-1}C^\top (C F_{t|t-1} C^\top +\Sigma_V)^{-1}.
\end{align}

We now consider the cost $J(K_{\mathsf{dyn}})$ (\ref{eq:costdefpo}) evaluated at the optimal filter (\ref{eq:filtereq}) and with variable $K$. With some abuse of notation, we denote this quantity $J(K;G)$ where $u_t$ is given by $u_t = K\hat x_t$, and $\hat x_t$ is defined by the Kalman filter (\ref{eq:filtereq}). We shall call the quantity $J(K;G)$ the restricted cost function, and note that it has almost the exact same form as the fully observed cost (\ref{eq:costdef}).  With these preliminaries established, the following lemma is straightforward to verify using Lemmas~\ref{lem:lqrcostlem} and \ref{lem:pgexplicitlemma} (and justifies the abuse of notation $J(K;G)$).

\begin{lemma}
\label{lem:pobsprels}
Consider a partially observed system $G=(A,B,C)$ of the form (\ref{eq:polds}). Then the restricted cost function satisfies
\begin{align*}
    J(K;G) = \tr P_K \Sigma_{\nu,G} +\tr Q (I-L_GC) 
\end{align*}
where $L_G$ and $\Sigma_{\nu,G}$ are the steady state quantities corresponding to recursions (\ref{eq:fric}) and (\ref{eq:fgain}) respectively and where $P_K$ as before is given by the Lyapunov equation (\ref{eq:lyapk}). Moreover, the policy gradient is given by
\begin{align}
\label{eq:pgpo}
   \nabla_{ K} J(K;G)= 2\left((R+B^\top P_K B)K+B^\top P_K A \right) \Gamma_{K,\nu,G}
\end{align}
where 
\begin{align}
\label{eq:gammadefpo}
    \Gamma_{K,\nu,G} \triangleq  \sum_{t=0}^\infty (A+BK)^t\Sigma_{\nu,G}(A+BK)^{t,\top}.
\end{align}
\end{lemma}

In other words, near the optimal controller $K_{\mathsf{dyn},_\star}$ the gradient with respect to the filter gain $K$ has the same form as in the state-feedback setting (\ref{eq:lds}). However, we stress at this point that neither the realization of the system (\ref{eq:polds})  nor the realization of the controller (\ref{eq:poobscont}) is unique. To remedy this, we will later verify in a scalar setting that our lower bounds are invariant under similarity transformation (see equation (\ref{eq:invtozero})).

\paragraph{Recovering Theorem~\ref{thm:mainthm}}
In the partially observed setting, we keep the exploration budget constraint (\ref{eq:budgeteq}) but the observation model is necessarily different. Namely, we assume that the learner instead has access to input-output data of the form $(y_{0,n},u_{0,n},\dots,u_{T-1,n},y_{T,n}), n \in [N]$.

In the partially observed setting, we thus define the analogous local minimax complexity as
\begin{equation}
    \begin{aligned}
\label{eq:minimaxdefpo}
    \mathfrak{M}_d(\e;G) 
    \triangleq \inf_{\reallywidehat{\nabla J} }\sup_{G': d(G,G')\leq \e} \E_{G'} \left\| \nabla_K J(K(G);G')-\reallywidehat{\nabla J} \right\|_{\mathsf{op}}
    \end{aligned}
\end{equation}
where the infimum is taken over all measurable functions of the data $(y_{0,n},u_{0,n},\dots,u_{T-1,n},y_{T,n}), n \in [N]$, $\nabla_K J(K;G)$ is given by equation (\ref{eq:pgpo}) and $d$ again is a metric on system parameters $G=(A,B,C)$. 

Equipped with the definition (\ref{eq:minimaxdefpo}) and Lemma~\ref{lem:pobsprels} the proof of the following result follows similarly to that of Theorem~\ref{thm:mainthm}. 
\begin{theorem}
\label{thm:mainthmpo}
Consider two systems $G_1 =(A,B,C)$ and $G_2(\Delta)=(A',B',C')$ with $A'=A-\Delta K_\star$, $B'=B+\Delta$ and $C'=C$. Then the local minimax complexity of estimating gradients (\ref{eq:minimaxdefpo}) is lower bounded as
\begin{equation}
\label{eq:mainthmpoeq}
\begin{aligned}
    \mathfrak{M}_d(\e; G_1)   
    \geq \sup_{d(G_1,G_2(\Delta))\leq \e} \left\| \Delta^\top P_{K_\star,G_1}(A+BK)  \Gamma_{K_\star,\nu,G_1}\right\|_{\mathsf{op}}
     \times \left(1-\sqrt{\frac{1}{2} d_{\mathsf{KL}}(G_1,G_2(\Delta))} \right).
\end{aligned}
\end{equation}
Above $d_{\mathsf{KL}}(G_1,G_2(\Delta))$ is the divergence between the induced probability measures over input-output data $(y_{0,n},u_{0,n},\dots,u_{T-1,n},y_{T,n}), n \in [N]$ between models $G_1$ and $G_2$.
\end{theorem}

By the data-processing inequality, the lower bound (\ref{eq:mainthmpoeq}) can be brought onto the exact same form as the lower bound (\ref{eq:mainthmpoeq}). Namely, we observe that\footnote{To see this, simply observe that $(y_0,\dots,y_T)$ is a stochastic function of $(x_{0,n},u_{0,n},\dots,u_{T-1,n},x_{T,n})$.}
\begin{align*}
    d_{\mathsf{KL}}(G_1,G_2(\Delta))  \leq  d_{\mathsf{KL}}(S_1,S_2(\Delta))
\end{align*}
where $ d_{\mathsf{KL}}(S_1,S_2(\Delta))$ is a slight overload of notation for the divergence between state-input data \\ $(x_{0,n},u_{0,n},\dots,u_{T-1,n},x_{T,n}), n \in [N]$ between models $G_1$ and $G_2$. In other words, all the results of Section~\ref{sec:mainsec} apply with $\Gamma_{K_\star,S_1}$ defined by equation (\ref{eq:gammadef}) exchanged for $\Gamma_{K_\star,\nu,G_1}$ defined in equation (\ref{eq:gammadefpo}) and $\Sigma_w$ exchanged for $\Sigma_t$ given by equation (\ref{eq:kalmancov}). While this is true for a fixed parametrization $G=(A,B,C)$, one may wonder whether the lower-bound relies on fundamental system-theoretic quantities or is simply a consequence of poor parametric choice for computing gradients. In the next example we show that the lower bound (\ref{eq:mainthmpoeq}) captures control-theoretic limitations that are independent of the state-space representation (\ref{eq:polds}). 

\paragraph{Bad Markov Parameters Imply Noisy Gradients}

Consider the almost scalar system $g_1=(a,B,c)$ given by
\begin{equation}
    \begin{aligned}
    \label{eq:scalarpobs}
    x_{t+1}&=ax_t + \begin{bmatrix} b & 0 \end{bmatrix} u_t +w_t\\
    y_t &= cx_t + v_t
    \end{aligned}
\end{equation}
defined consinstently with system (\ref{eq:polds}), but specifically $a,b,c \in \mathbb{R}$ and $Q=\Sigma_V=\Sigma_W =1$ and $R = I_2$. Note that the maximum singular values of the first Markov parameter of $g_1$ is equal to the product $m=cb$ and that this is $m$ is invariant under similarity transformation. We consider the two systems $g_1=(a,B,c)$ and $g_2=(a,B(\Delta),c)$ and where $B(\Delta)= \begin{bmatrix} b & \Delta \end{bmatrix}$. Observe that the optimal policy to system (\ref{eq:scalarpobs}) is of the form $K_\star = \begin{bmatrix} k_\star & 0 \end{bmatrix}$ and that the gramians $P_{K_\star,\nu,g_1}$ and $\Gamma_{K_\star,g_1}$ are scalar and equal to $P_{K_\star,\nu,g_1}=P_{k_\star,\nu,g_1}$ and $\Gamma_{K_\star,g_1}=\Gamma_{k_\star,g_1}$ respectively. In other words, the second input has no effect on the system (\ref{eq:scalarpobs}), but as well shall see, $g_2$ is very sensitive to perturbations $\Delta$ whenever the largest singular value of the Markov parameter $m=|cb|$ is small.

If we denote by $d_{\mathsf{KL}}(s_1,s_2)$ the KL divergence between scalar input-state trajectories drawn from $g_1$ and $g_2$, we have by Lemma~\ref{lem:KLcalc} that
\begin{align*}
    d_{\mathsf{KL}}(s_1,s_2) &= \sum_{n=1}^N \sum_{t=0}^{T-1}\E_{g_1} \frac{1}{2}( \Delta u_t)^2&& \textnormal{(Lemma~\ref{lem:KLcalc})} \\
    &\leq   \frac{1}{2}\Delta^2 \sum_{n=1}^N \sum_{t=0}^{T-1}\E_{g_1}  u_t^2\\
    &\leq  \frac{1}{2}\Delta^2 NT\beta && \textnormal{(by (\ref{eq:budgeteq}))} \\
    &\leq \frac{1}{2} &&\left(\textnormal{ if } \Delta^2 = \frac{1}{ NT \beta }\right).
\end{align*}
Invoking Theorem~\ref{thm:mainthmpo} this implies the local minimax lower bound
\begin{equation}
\begin{aligned}
\label{eq:damnedifdo}
    \mathfrak{M}_d(\e_{N,T}; g_1) &\geq \frac{1}{2\sqrt{ NT \beta} } P_{k_\star,\nu,g_1}(a+bk_\star) \Gamma_{k_\star,\nu, g_1} 
\end{aligned}
\end{equation}
where $\e_{N,T} \asymp 1/\sqrt{NT}$. Inequality (\ref{eq:damnedifdo}) in itself is an instance specific lower bound for scalar partially observed systems of the form (\ref{eq:scalarpobs}). Further, the inequality implies that if the Markov parameter $|m|=|cb|$ is small, estimating gradients is always hard. Namely, we make the following observations\footnote{To verify these claims, observe that the scalar quantities $P_{k_\star,\nu,g_1}$, $k_\star$ and $\Sigma_{k_\star,\nu,g_1}$ have closed form solutions.}:
\begin{itemize}
    \item $P_{k_\star,\nu,g_1}$ tends to infinity at rate $1/b^2$ as $b$ tends to $0$. Moreover, $P_{k_\star,g1}$ is always lower-bounded by $1$.
    \item The large and small $c$ asymptotics of $\Sigma_{k_\star,\nu,g_1}$ are proportional to $1/c^2$
    \item The factor $|a+bk_\star|$ tends to 0 no faster than $1/b^2$ and tends to $\min(1,|a|)$ as $b\to 0$ (and $|a|$ is invariant under similarity transform).
\end{itemize}
Combining these observations, we see that as the system invariant $|m|=|cb|$ tends to zero, gradients become arbitrarily noisy; the lower bound (\ref{eq:damnedifdo}) tends to infinity. In other words, we have established that \begin{align}
\label{eq:invtozero}
    \liminf_{|cb|\to 0} \mathfrak{M}_d(\e_{N,T}; g_1) = \infty.
\end{align}
Thus, we obtain an RL analogue to the well-known fact that reparametrization cannot help controlling a partially observed system as any gain in observability is offset by a proportional loss in controllability and vice versa.

\section{Discussion}

In this work we showed that estimating policy gradients can become arbitrarily hard due to known control-theoretic fundamental limitations \cite{doyle1978guaranteed} by leveraging the classic two point method due to Le Cam \cite{lecam1973convergence}. For instance, we showed with system (\ref{eq:scalarpobs}) that a partially observed system with small Markov parameters necessarily has noisy policy gradients and that this holds independently of the parametrization. Our bounds also show that learning controllers that are close to marginal stability can be hard. This is similar to what has already been observed for adaptive LQR/LQG in \cite{ziemann2022regret}. Leveraging results from \cite{Tasosunpub} we further show that estimating policy gradients can suffer from exponential complexity in the system dimension. From a broader perspective, these results work toward elucidating when learning to control is feasible.

\textit{Acknowledgements:} 
Ingvar Ziemann and Henrik Sandberg are supported by the Swedish Research Council (grant 2016-00861) and the Swedish Foundation for Strategic Research through the CLAS project (grant RIT17-0046). Nikolai Matni is supported in part by NSF awards CPS-2038873 and CAREER award ECCS-2045834, and a Google Research Scholar award.


\bibliographystyle{unsrtnat}
\bibliography{main}

\appendix
\section{Auxilliary Lemmas}

\paragraph{Proofs for the Preliminary Results}

\textit{Proof of Lemma~\ref{lem:lqrcostlem}.}
For $T\in \mathbb{N}$ and law $u_t = Kx_t$ we have
\begin{equation}
\begin{aligned}
    \sum_{t=0}^{T-1} x_t^\top Q x_t + u_t^\top R u_t &=\sum_{t=0}^{T-1} x_t^\top (Q+K^\top RK) x_t\\
    &= \sum_{t=0}^{T-1} \tr \left[ (Q+K^\top RK) x_tx_t^\top \right].
\end{aligned}
\end{equation}
Recall $\Gamma_K = \sum_{t=0}^\infty (A+BK)^t\Sigma_W(A+BK)^{t,\top}$ so that $
    \E x_t x_t^\top = \Gamma_K +o(1)$, where $o(1)$ tends to $0$ as $t$ tends to infinity. Hence, by averaging and taking limits we see that
\begin{equation}
\begin{aligned}
    &J(K)\\ 
    &= \tr \left[(Q+K^\top R K) \Gamma_K \right]\\
    &= \tr \left[(Q+K^\top R K)  \sum_{t=0}^\infty (A+BK)^t\Sigma_W(A+BK)^{t,\top} \right]\\
    &=\tr \left(\sum_{t=0}^\infty \left[(A+BK)^{t,\top}(Q+K^\top R K)   (A+BK)^t \right]\Sigma_W\right).
\end{aligned}
\end{equation}
The result follows since $(A+BK)$ is stable by hypothesis.
\hfill $\blacksquare$

\textit{Proof of Lemma~\ref{lem:pgexplicitlemma}.} 
Fix a controller $K$. By Lemma~\ref{lem:lqrcostlem} the average cost of the system (\ref{eq:lds}), $J(K)$, is the same as the total cost of the deterministic system $\tilde x_{t+1} = A\tilde  x_t +B \tilde u_t$ with $\tilde x_0 \sim N(0,\Sigma_W)$:
\begin{align*}
    J_S(K) = \E_{K,S} \sum_{t=0}^\infty \tilde x_t^\top Q \tilde x_t + \tilde u_t R u_t.
\end{align*}The result now follows by Lemma 1 of \cite{fazel2018global}.\hfill $\blacksquare$

\paragraph{Information-Theoretic Lower Bounds}

Intuitively, estimating a function $f(A,B)$ while only having access to samples from the unknown system (\ref{eq:lds}) becomes hard if there are parameter variations $A'$ and $B'$ such that the behavior of $x'_{t+1} =A'x'_t+B'u'_t+w_t$ is very close to that of system (\ref{eq:lds}) while the difference between $f(A,B)$ and $f(A',B')$ is large. This can be formalized by the Le Cam's two-point method:

 \begin{lemma}[Le Cam's Two Point Method]
\label{thm:lecam}
Fix two sets $\mathcal{M} $ and $ \mathcal{D} $. Let $L: \mathcal{M}\times  \mathcal{D} \to \mathbb{R}_+$ be any loss function and suppose that $S_1,S_2 \in  \mathcal{M}$ satisfy $L(S_1,\mathsf{dec})+L(S_2,\mathsf{dec}) \geq \delta, \: \forall \mathsf{dec} \in \mathcal{D}$. Then
 \begin{align*}
     \inf_K \sup_{S\in \mathcal{S}} \E_S L(S,K) \geq \frac{\delta}{2} (1-d_{\mathsf{TV}}(P_{S_1},P_{S_2})).
 \end{align*}
 \end{lemma}
 
 In other words, if for any decision the average loss is large, then a decision-maker that cannot distinguish between these two instances will suffer large loss on average, and therefore also in the worst case.

\textit{Proof of Lemma~\ref{thm:lecam}.}
 We lower-bound the supremum over $\mathcal{M}$ by an expectation over the two-point mixture distribution supported on $S_1$ and $S_2$ as follows:
 \begin{equation*}
     \begin{aligned}
     \inf_{\mathsf{dec}} \sup_{S\in \mathcal{\bar S}} \E_S L(S,\mathsf{dec}) &\geq \inf_{\mathsf{dec}}  \frac{\E_S L(S_{1},\mathsf{dec})+\E_{S_2} L(S_2,\mathsf{dec})}{2} \\
     &=\frac{1}{2}\inf_{\mathsf{dec}}\Bigg(\int L(S_1,\mathsf{dec}(x))P_{S_1}(dx) \\
     &+\int L(S_2,\mathsf{dec}(x))P_{S_2}(dx)  \Bigg) \\
    &\geq \frac{\delta }{2}\left(\int \min(P_{S_1}(dx), P_{S_2}(dx))  \right)\\
    &=\frac{\delta}{2} (1-d_{\mathsf{TV}}(P_{S_1},P_{S_2}))
    \end{aligned}
 \end{equation*}
 as per requirement.\hfill $\blacksquare$

Since $d_{\mathsf{TV}} \leq \sqrt{\frac{1}{2} d_{\mathsf{KL}}}$ by Pinsker's inequality, the following result is convenient to state.

\begin{lemma}
\label{lem:KLcalc}
Let $S_0 = (A_0,B_0)$ and $S_1=(A_1,B_1)$ and denote by $P_1$ and $P_2$ the induced probability measures over samples $(x_0,\dots,x_T)$ satisfying the recursion (\ref{eq:lds}) with $A=A_i,B=B_i, i =0,1$. Then
\begin{align*}
    d_{\mathsf{KL}}(P_1,P_2) = \sum_{t=0}^{T-1} \frac{1}{2} \E_1\|(A_0-A_1)x_t+(B_0-B_1)u_t\|_{\Sigma_W^{-1}}^2
\end{align*}
where $\E_1$ denotes integration with respect to $\mathbf{P}_1$ and the norm $\|\cdot\|_{\Sigma_W^{-1}}$ is the Mahalanobis norm with kernel $\Sigma_W^{-1}$.
\end{lemma}

\textit{Proof of Lemma~\ref{lem:KLcalc}.}
Each random variable $x_t$ is conditionally Gaussian given $(x_0,\dots,x_{t-1})$ so that the KL divergence is given by half the expected square difference in conditional mean. The result follows by straighforward computation and the chain rule. See for instance \cite[Chapter 8]{cover1999elements}.
\hfill $\blacksquare$

\begin{figure}
\begin{subfigure}{.5\textwidth}
  \centering
  \includegraphics[width=.8\linewidth]{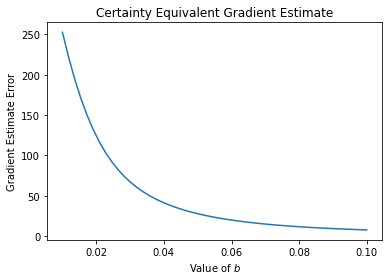}  
\end{subfigure}
\begin{subfigure}{.5\textwidth}
  \centering
  \includegraphics[width=.8\linewidth]{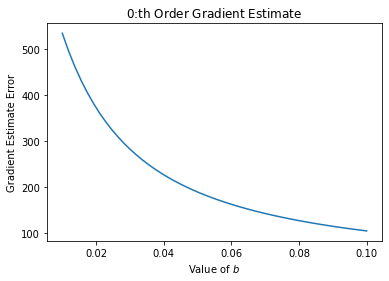}  
\end{subfigure}
\caption{Gradient estimate spread as a function of $b$ for the scalar system~(\ref{eq:sysscalar}). Notice that poor controllability (small $b$), leads to noisy gradients. The vertical axes show the standard deviation of $\left\| \nabla_K J(K;S)-  \widehat{\nabla_K J} \right\|_{\mathsf{op}}$ across multiple trajectories.}
\label{fig:plots}
\end{figure}

\paragraph{Simulation}
In Figure~\ref{fig:plots} we numerically verify the claims made for scalar systems in Section~\ref{subsec:cons} with $a=1$, $\sigma_w=1$ and variable $b$. For the first plot, we use trajectories of length $T=100$ and compute the least squares certainty equivalent (plug-in) gradient estimate using a single trajectory. The error is then averaged over $N=100$ trajectories. For the second plot, we also use trajectories of length $T=100$. However, we use $N=10000$ many trajectories divided into batches, with each batch containing $100$ trajectories. Each batch is then used to compute a $0$-th order gradient estimate (see \cite[Algorithm 1]{fazel2018global}). The second plot shows the estimator error averaged over these batches. Notice that for either estimator, the performance diverges in the low-controllability regime $b \approx 0$.

\end{document}